\documentclass{amsart}

\usepackage[dvips]{graphicx}

\newtheorem{theorem}{Theorem}[section]
\newtheorem{lemma}[theorem]{Lemma}
\newtheorem{corollary}[theorem]{Corollary}
\newtheorem{proposition}[theorem]{Proposition}

\theoremstyle{definition}
\newtheorem{definition}[theorem]{Definition}
\newtheorem{example}[theorem]{Example}
\newtheorem{question}[theorem]{Question}

\newcommand{\dt}{\mathrm{DT}}

\begin{document}

\title[Complexes of Directed Trees]
  {Complexes of Directed Trees and
  Independence Complexes}

\author
  {Alexander Engstr\"om}

\address
  {Department of Computer Science, 
  Eidgen\"ossische Technische Hochschule, 
  Z\"urich, Switzerland.} 

\email
  {engstroa@inf.ethz.ch}


\date\today


\thanks
  {Research supported by ETH and Swiss National 
  Science Foundation Grant PP002-102738/1}

\begin{abstract}
  The theory of complexes of directed trees
  was initiated by Kozlov to answer a
  question by Stanley, and later on, results 
  from the
  theory were used by Babson and Kozlov in their
  proof of the Lov\'asz conjecture. We develop
  the theory and prove that complexes on 
  directed acyclic graphs are shellable.

  A related concept is that of independence
  complexes: construct a simplicial complex on
  the vertex set of a graph, by including 
  each independent set of vertices as a simplex.
  Two theorems used for breaking and gluing
  such complexes are proved and applied to 
  generalize results by Kozlov.

  A fruitful restriction is
  anti-Rips complexes: a subset $P$ of a 
  metric space is the vertex set of the complex,
  and we
  include as a simplex each subset of $P$ with no
  pair of points within distance $r$.
  For any finite subset $P$ of $\mathbb{R}$ the
  homotopy type of the anti-Rips 
  complex is determined.
\end{abstract}

\maketitle

\section{Introduction}
  The theory of complexes of directed trees was 
  initiated by Kozlov \cite{K1} and Babson and
  Kozlov used results from this theory in their 
  proof of the Lov\'asz conjecture \cite{BK}. 
  In this paper we study three
  ways to connect topology with combinatorics
  by constructing simplicial complexes: 
\begin{center}
  \begin{tabular}{l|l|l}
  Name & Base object & Restriction on simplices\\
  \hline
  $\dt(G)$ & Directed graph $G$ & They are 
  directed forests of $G.$\\
  $\mathrm{Ind}(G)$ & Undirected graph $G$ & No
  vertices are adjacent in $G$.\\
  $\mathrm{AR}_r(P)$ & Point set $P$ in metric
  space & No two vertices within distance $r$.
  \end{tabular}
\end{center}
  To decide the homotopy type of complexes on
  directed graphs, Kozlov \cite{K1} constructed
  complexes on undirected graphs.
  In the second part of
  this paper we show that by moving the
  problems into complexes on point sets of 
  metric spaces, the homotopy type can be
  determined for several classes, and
  naturaly generalized. 
  The first section treats
  complexes on directed trees. It is shown
  that certain complexes are shellable,
  for example those on directed acyclic graphs.
\section{Complexes of directed trees}
  In this section all graphs are directed.
  A \emph{directed 
  forest} is an acyclic graph with at most one 
  edge directed to each vertex.
  Equivalently,
  a directed forest is a collection of disjoint
  directed trees with all edges oriented away
  from the root.
\begin{definition}
  Let $G$ be a directed graph. The complex of
  directed trees, $\dt(G)$ have the edges of
  $G$ as vertex set. A set of edges is a 
  simplex of $\dt(G)$ if the edges viewed as a
  graph is a directed forest.
\end{definition}
  Constructing simplicial complexes from graphs
  can, in principal, be done in two ways: For
  graphs that satisfy certain properties
  either their edges or their vertices 
  form a simplex. For complexes of directed trees
  it is the edges, but later on complexes of
  the other kind are used.

  A directed forest $H\subseteq G$ is 
  \emph{maximal} if $H'$ is not a directed forest
  for any $H\subset H' \subseteq G$.
  The \emph{roots} of a directed forest are the 
  roots of the trees in the forest. A maximal 
  face of $\dt(G)$ is the edge set of a 
  maximal directed forest in $G$, and the other 
  way around.
\begin{definition} Let $R\subseteq V(G)$. The 
  simplicial complex $\dt_R(G)\subseteq \dt(G)$ 
  is generated by the faces of $\dt(G)$ which are
  edge sets of  directed forests with $R$ as 
  roots. 
\end{definition}
  The forests in $\dt(G)$ with roots $R$ are 
  the maximal faces of $\dt_R(G)$. 
\begin{definition}
  An edge $(x\rightarrow y)$ of $G$ is 
  \emph{nice} in a subcomplex $\Delta$ of
  $\dt(G)$ if
\begin{itemize}
  \item[\emph{(i)}]
  there is an edge $(z\rightarrow y)$ in 
  $\Delta$ such that $z\neq x$;
\item[\emph{(ii)}]
  any forest $F\in\Delta$ without an edge
  directed to $y$ can be extended with
  $(x\rightarrow y)$, and 
  $F\cup\{(x\rightarrow y)\}\in\Delta$.
\end{itemize}
\end{definition}

\begin{example}\label{example:1}
  Let $G$ be the directed graph in Figure
  \ref{fig:pix1}.
\begin{figure}
  \begin{center}
  \includegraphics*{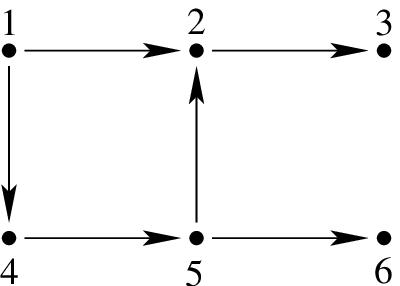}
  \end{center}
  \caption{Graph $G$.}
  \label{fig:pix1}
\end{figure}
  Let us find some nice edges in $\dt_R(G)$, 
  where $R=\{1,4\}$. The maximal faces are
  drawn in Figure \ref{fig:pix23}.
\begin{figure}
  \begin{center}
  \begin{tabular}{ccc}
  \includegraphics*{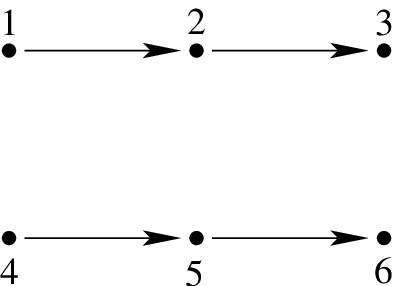} & $\quad$ &
  \includegraphics*{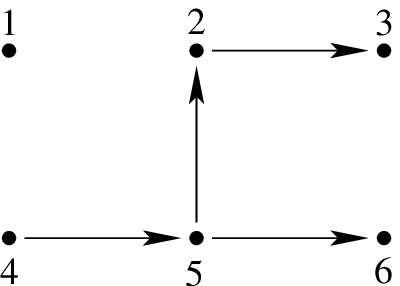} \\
  \end{tabular}
  \end{center}
  \caption{Maximal faces $F_1$ and $F_2$ of
  $\dt_R(G)$.}
  \label{fig:pix23}
\end{figure}
  We partition the vertices into left and right 
  side of a dotted line: A vertex $v$ is on the 
  left side if the tree which $v$ is in have the
  same root $r$ in all maximal faces of 
  $\dt_R(G)$, and the
  path from $r$ to $v$ is the same in all
  maximal faces of $\dt_R(G)$.
  All other vertices are on the right side. The
  union of $F_1$ and $F_2$ with their vertices
  partitioned is depicted in Figure 
  \ref{fig:pix4}.
\begin{figure}
  \begin{center}
  \includegraphics*{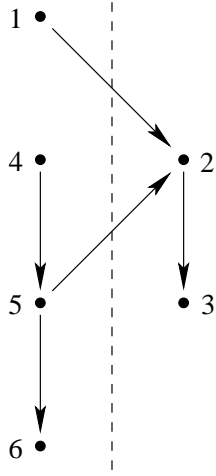}
  \end{center}
  \caption{The union of $F_1$ and $F_2$ with
  their vertices partitioned.}
  \label{fig:pix4}
\end{figure}
  In this example it is not hard to see that the 
  edges crossing the dotted line are nice. If 
  all maximal faces have the same roots, this is 
  true in general.
\end{example}
\begin{proposition}\label{prop:1}
  All edges crossing the dotted line in a
  construction as in Example \ref{example:1}
  are nice if the maximal faces have the same
  roots.
\end{proposition}
\begin{proof}
  First we need to show that an edge 
  $(z\rightarrow w)$ can only cross the dotted 
  line from left to right. Assume the contrary,
  i.e. that 
  $(z\rightarrow w)$ is from right to left. 
  Since $w$ is on the left, the tree which $w$ 
  is in has the same root $r$ in all maximal 
  faces, 
  and the path from $r$ to $w$ is the same.
  The vertex before $w$ in that path is always 
  the same, and $(z \rightarrow w)$ is in a 
  maximal face, so the vertex is $z$. But then 
  the path from $r$ to $z$ is the same in 
  all maximal faces, and $z$ should be on the 
  left side, which is a contradiction.

  Assume that $(x\rightarrow y)$ is an edge 
  crossing the dotted line and that the maximal 
  faces are $\{ F_i \}_{i\in I}$.

  The vertex $y$ is not a root, hence there is 
  an edge to $y$ in all maximal faces. If 
  $(x\rightarrow y)$ was the only edge
  to $y$ in
  $\cup_{i\in I} F_i$, it would be in all 
  maximal faces. But then the path to $y$ is the 
  same in all maximal faces, and $y$ would
  be on the 
  left side. Since $y$ is on the right side,
  there is an edge $(z\rightarrow y)$ such
  $z\neq x$. Therefore condition (i) of
  the definition of nice edges is satisfied.

  If $(x\rightarrow y)$ does not induce a 
  directed cycle when added to a face, then 
  condition (ii) is fulfilled. A directed cycle 
  which contains $(x\rightarrow y)$ as an edge
  would cross the dotted line at least twice. 
  But all edges crossing 
  the dotted line go from left to right, hence 
  $(x\rightarrow y)$ cannot induce a directed 
  cycle.
\end{proof}

\begin{lemma}\label{lemma:findNice}
  If $R\subseteq V(G)$ is nonempty, and 
  $\dt_R(G)$ has more than one maximal face, 
  then there is an edge $(x\rightarrow y) 
  \in E(G)$ which is nice in $\dt_R(G)$.
\end{lemma}
\begin{proof}
  Construct the left/right-partition of the 
  vertices as in Example \ref{example:1}.
  If there are no edges crossing the dotted line,
  all edges are on the left side, and
  there is only one maximal face.
  Hence there is an edge 
  crossing the dotted line, and by Proposition
  \ref{prop:1} that edge is nice.
\end{proof}
\begin{definition}
  A simplicial complex $\Delta$ is 
  \emph{shellable} if 
  its maximal faces can be ordered $F_1,F_2, 
  \ldots ,F_n$ such that for all 
  $1\leq i<k\leq n$, there are $1\leq j<k$ and 
  $e\in F_k$, such that 
  $F_i\cap F_k \subseteq F_j
  \cap F_k = F_k\setminus \{e\}$.
\end{definition}

\begin{lemma}\label{lemma:constructMax}
  Let $F_i$ and $F_k$ be maximal faces of
  $\dt_R(G)$, and $(x\rightarrow y)\in
  F_i \setminus F_k$ a nice edge in
  $\dt_R(G)$. Then there is a maximal face
  $F_j$ of $\dt_R(G)$, and an  
  $e\in F_k$, such that
  $F_i\cap F_k \subseteq F_j
  \cap F_k = F_k\setminus \{e\}$.
\end{lemma}
\begin{proof}
  There is an edge $(z\rightarrow y)$ in $F_k$,
  since it is maximal and $y\not\in R$.
  Replace it with $(x\rightarrow y)$ to
  construct $F_j$. From the definition of
  nice edges we have that $F_j\in\dt_R(G)$.
  Now $F_j \cap F_k = F_k \setminus 
  \{(z\rightarrow y)\}$. Since $(x \rightarrow 
  y)\in F_i$ and $(z \rightarrow 
  y)\not\in F_i$, we conclude that 
  $F_i\cap F_k = F_i \cap (F_k \setminus 
  \{(z\rightarrow y)\}) \subseteq $ $F_k 
  \setminus \{(z\rightarrow y)\}$ $= F_j\cap 
  F_k$.
\end{proof}

\begin{theorem}\label{theorem:rootShellable} 
  $\dt_R(G)$ is shellable.
\end{theorem}
\begin{proof}
  The proof is by induction over the number of 
  maximal faces. If $\dt_R(G)$ has one maximal 
  face, it is a simplex, and thus shellable.
  Assume that $\dt_R(G)$ has more that one 
  maximal face. By Lemma~\ref{lemma:findNice}
  there is a nice edge $(x \rightarrow y)$. 
  Define $G',G''\subset G$ as follows:
  \[ E(G')=E(G)\setminus \{(z \rightarrow y)\mid
  z\neq x\} \quad E(G'')=E(G)\setminus 
  \{(x \rightarrow y)\}.\]
  A maximal face of $\dt_R(G)$ has exactly one 
  edge to $y$, thus the set of maximal faces of 
  $\dt_R(G)$ is the disjoint union of the sets 
  of maximal faces of $\dt_R(G')$ and $\dt_R
  (G'')$. Since $(x \rightarrow y)$ is nice, it
  is in some, but not in all maximal faces of
  $\dt_R(G)$. Both $\dt_R(G')$ and $\dt_R(G'')$
  have smaller numbers of maximal faces than 
  $\dt_R (G)$, so by induction they are 
  shellable. 
  Order the maximal faces of $\dt_R
  (G):$
  \[ \underbrace{F_1, F_2, \ldots, F_t,}
  _{\text{Shelling order of }\dt_R(G'),\:\:
   (x \rightarrow y)\in F_l.\:\:}
  \underbrace{ F_{t+1}, F_{t+2}, \ldots, F_{t+s}}
  _{\text{  Shelling order of }\dt_R(G''),
  \:\:(x \rightarrow y)\notin F_l.}\]
  If for all $1\leq i<k\leq s+t$, there are 
  $1\leq j < k$ and $(z\rightarrow w)\in F_k$ 
  such that $F_i \cap F_k \subseteq F_j\cap F_k =
   F_k\setminus \{ (z\rightarrow w) \}$, then 
  $\dt_R(G)$ is shellable.
  If $1\leq i <k 
  \leq t$ or $t<i<k \leq s+t$ we are done. 
  Assume $1\leq i \leq t < k \leq s+t$. 

  The nice edge $(x \rightarrow y)$ is in $F_i$, 
  but not in $F_k$. Construct $F_j$ as described 
  in Lemma~\ref{lemma:constructMax}. 
  The edge $(x \rightarrow y)$ is in $F_j$, 
  so $j\leq t<k$, and we have a shelling order.
\end{proof}

  A \emph{directed acyclic graph} is a directed 
  graph without directed cycles. It is not hard 
  to see that the vertices of a directed acyclic 
  graph $G$ can be ordered so that if $(x
  \rightarrow y)\in E(G)$, then $x$ is before 
  $y$.
  This is usually called the topological order. 
  Let $R$ be the set of vertices of $G$ with 
  no edges directed to them.
  The first vertex in the order is in $R$, 
  so the set is not empty. Since $G$ is a 
  directed acyclic graph, so are all its 
  subgraphs.
  Hence the maximal subgraphs, with at most one 
  edge to each vertex, are the maximal forests. 
  All maximal forests contain an edge to each 
  vertex in $V(G) \setminus R$. 
  Hence the roots of all 
  maximal forests are $R$, and 
  $\dt(G)=\dt_R(G)$.
\begin{corollary}\label{corollary:shellableDAG}
  If $G$ is a directed acyclic graph, then $\dt
  (G)$ is shellable. 
\end{corollary}

  A pure shellable simplicial 
  complex $\Delta$ is homotopy 
  equivalent to a wedge of spheres of the same 
  dimension as $\Delta$, or it is contractible. 
  Thus, by 
  calculating the Euler characteristic of 
  $\Delta$, its homotopy type can be determined.
  See Bj\"orner and Wachs, \cite{BW}, for a 
  proof of that, 
  and further extensions.

  Denote with $d^-(v)$ the number of edges 
  directed to a vertex $v$.
\begin{lemma}\label{lemma:eulerDAG}
  If $G$ is a directed acyclic graph with at
  least one edge, then 
  \[\tilde{\chi} (\dt(G)) 
  =-\prod_{v\in V(G)\setminus R} (1-d^-(v)),\]
  where $R$ is the set of vertices without edges
  directed to them.
\end{lemma}
\begin{proof}
  The proof is by induction over the number of
  vertices not in $R$. If there is only one vertex
  not in $R$, then the
  complex $\dt(G)$ is homotopy equivalent to
  $d^-(v)$ disjoint points, and the formula is 
  true.
 
  If $V(G)\setminus R$ has more than two 
  vertices,  
  order the vertices so that if $(x\rightarrow y)
  \in E(G)$, then $x$ is before $y$, and so that 
  all vertices of $R$ come before the other 
  ones. Denote the last vertex 
  in this order with $w$. Let $G'$ be the 
  induced subgraph of $G$ with vertex set 
  $V(G)\setminus \{ w\}$. Since $w$ is the last 
  one ordered, there are no edges from 
  $w$, and $E(G)\setminus E(G')$ are the 
  $d^-(w)$ edges to $w$.

  Let $\alpha (i)$ be the number of subgraphs of 
  $G$ which are forests with $i$ edges, and 
  $\alpha '(i)$ similarly for $G'$. A forest 
  with $i$ edges in $G$ either has no edge to 
  $w$, or one of the $d^-(w)$ edges to $w$, 
  hence for $i>0$
  \[ \alpha(i)=\alpha '(i)+d^-(w)\alpha'(i-1).\]
  Clearly $\alpha (0)=\alpha'(0)=1$. 
  The reduced 
  Euler characteristic of $\dt (G)$ is
  \[ \begin{array}{rcl}
  \tilde{\chi} (\dt (G)) & = & \displaystyle
  \sum_{i\geq 0} (-1)^{i+1} \alpha(i) \\
  &=& -1 +\displaystyle \sum_{i\geq 1}(-1)^{i+1} 
  \left(\alpha'(i) + d^-(w) \alpha'(i-1)\right)\\
  & = & -1+\displaystyle \sum_{i\geq 1}(-1)^{i+1}
  \alpha'(i)+ \displaystyle  d^-(w)\sum_{i\geq 1}
  (-1)^{i+1} \alpha '(i-1) \\
  &=& \displaystyle \sum_{i\geq 0}(-1)^{i+1}
  \alpha'(i)-\displaystyle  d^-(w)\sum_{i\geq 0}
  (-1)^{i+1}  \alpha '(i) \\
  &=& (1- d^-(w)) \displaystyle \sum_{i\geq 0}
  (-1)^{i+1}  \alpha '(i) \\
  &=& (1- d^-(w)) \tilde{\chi} (\dt (G'))
  \end{array} \]
  Substituting
  the formula for $ \tilde{\chi} 
  (\dt(G'))$ concludes the proof.
\end{proof}

\begin{theorem}\label{theorem:homEqDAG}
  If $G$ is a directed acyclic graph, then $\dt
  (G)$ is homotopy equivalent to a wedge of 
  $\prod_{v\in V(G)\setminus R} (d^-(v)-1)$ 
  spheres of dimension $\#V(G)-\#R-1$, where $R$ 
  is the set of vertices without edges
  directed to them.
\end{theorem}
\begin{proof}
  The complex $\dt(G)$ is shellable by 
  Corollary~\ref{corollary:shellableDAG}, and the
  reduced Euler characteristic is $\pm \prod_{
  v\in V(G) \setminus R }(d^-(v)-1)$ by 
  Lemma~\ref{lemma:eulerDAG}. The maximal faces 
  of $\dt(G)$ are the maximal forests of $G$. 
  Since the forests have edges exactly to the 
  vertices with non-zero in-degree, there are 
  $\#V(G)-\#R$ 
  edges in a maximal forest, and the 
  dimension of a maximal face is $\#V(G)-\#R-1$.
\end{proof}
\begin{corollary}
  If $G$ is a directed acyclic graph, then 
  the following statements are equivalent:
\begin{itemize}
  \item $\dt(G)$ is a cone with one of the edges
  of $G$ as apex.
  \item There is an edge of $G$ which is in all
  maximal forests.
  \item There is a vertex in $G$ with in-degree 
  1.
  \item The product of all $(d^-(v)-1)$ for 
  $v\in V(G)$ is zero.
  \item $\dt(G)$ is contractible.
\end{itemize}
\end{corollary}
\begin{proof}
  Follows directly from Theorem 
  \ref{theorem:homEqDAG} and its proof.
\end{proof}

\section{Independence complexes}
  In the previous section the graphs were 
  directed, but from now on all graphs are 
  assumed to be undirected.
  A subset of the vertex set of a
  graph is \emph{independent} if no two vertices
  in it are adjacent. In a graph $G$, the 
  \emph{neighborhood} of a vertex $v$, $N_G(v)$ 
  is the set of vertices which are adjacent to
  $v$.
  If it is clear which $G$ is meant, we just
  write $N(v)$.

  If $W\subseteq V(G)$ then $G[W]$ is the induced
  subgraph with vertex set $W$, and $G\setminus W
  =G[V(G)\setminus W]$. Similarly for simplicial
  complexes, if $W\subseteq \Delta^{(0)}$ then
  $\Delta[W]$ is the induced subcomplex, and 
  $\Delta\setminus W = \Delta[\Delta^{(0)}
  \setminus W]$.

\begin{definition}
  Let $G$ be an undirected graph. The 
  \emph{independence complex} of $G$, denoted 
  $\mathrm{Ind}(G)$, is a simplicial complex with
  vertex set $V(G)$, and $\sigma\in\mathrm{Ind}
  (G)$ if $\sigma$ is an independent set of $G$.
\end{definition}
  
  Proving the following standard facts is a
  good excercise to get acquainted with
  independence complexes.
\begin{itemize}
  \item If $A\subseteq V(G)$ then $\mathrm{Ind}
  (G[A]) = \mathrm{Ind}(G)[A]$.
  \item If $A,B\subseteq V(G)$ then
  $\mathrm{Ind}(G[A]) \cap \mathrm{Ind}(G[B]) = 
  \mathrm{Ind}(G[A\cap B])$.
  \item If $v \in V(G)$ then
  $\mathrm{lk}_{\mathrm{Ind}(G)}(v)=
  \mathrm{Ind}(G\setminus (N(v)\cup\{v\}))$, 
  and\\ 
  $\mathrm{st}_{\mathrm{Ind}(G)}(v)=
  \mathrm{Ind}(G\setminus N(v))=
  v\ast \mathrm{lk}_{\mathrm{Ind}(G)}(v)$.
  \item If $\mathrm{Ind}(G\setminus 
  (N(v)\cup\{v\}))$ is contractible, then
  $\mathrm{Ind}(G)\simeq \mathrm{Ind}(G
  \setminus v)$
  \item If $v \in V(G)$ then
  $\mathrm{Ind}(G)$ is the union of 
  $\mathrm{st}_{\mathrm{Ind}(G)}(v)$
  and $\cup_{w\in N(v)}
  \mathrm{st}_{\mathrm{Ind}(G)}(w)$.
  \item 
  If $v,w\in V(G)$ are adjacent and 
  $N(v) \setminus
  \{w\} \subseteq N(w)\setminus \{v\}$, then 
  $\mathrm{st}_{\mathrm{Ind}(G)}(v)
  \supseteq \mathrm{lk}_{\mathrm{Ind}(G)}(w)$.
\end{itemize}

  If the neighborhood of a vertex $v$ is
  included in the neighborhood of another vertex,
  the removal of $v$ from the graph is called 
  a fold. In the theory of ${\tt Hom}$--complexes
  folds is a fundamental tool for reducing the
  size of the input graphs while preserving the
  simple homotopy type, see \cite{K3}.
  But in 
  contrast with folds for ${\tt Hom}$--complexes,
  the vertex with the larger neighborhood is 
  removed to preserve 
  the simple homotopy type of an
  independence complex.
\begin{lemma}\label{lemma:fold}
  If $N(v)\subseteq N(w)$ then 
  $\mathrm{Ind}
  (G)$ collapses onto $\mathrm{Ind}(G\setminus 
  \{w\})$.
\end{lemma}
\begin{proof}
  Match each maximal $\sigma$, such that $w\in 
  \sigma$ and $v\not \in \sigma$, with 
  $\sigma \cup \{v\}$, and remove them by an
  elementary collapse step. Repeat this until
  all $\sigma$ such that $w\in \sigma$ are
  gone.
\end{proof}
  In particular, if $N(u)=\{v\}$ and $w\in N(v)
  \setminus\{u\}$, then $\mathrm{Ind}(G)$ 
  collapses onto $\mathrm{Ind}(G\setminus\{w\})$.
\begin{proposition}[\cite{EH}]
  If $G$ is a forest then $\mathrm{Ind}(G)$
  is either contractible or homotopy equivalent
  to a sphere.
\end{proposition}
\begin{proof}
  It is sufficient to show that successive use
  of Lemma \ref{lemma:fold} starting with $G$ 
  gives a graph
  $H$ with no adjacent edges, since  
  $\mathrm{Ind}(H)$ is either contractible or 
  homotopy equivalent to a sphere, and each use
  of the lemma provides a collapse.

  The proof is by induction on the number of 
  edges. If there are no adjacent edges we are
  done. If $G$ is a forest with some adjacent 
  edges, then there is a vertex $u$ with only 
  one neighbor $v$, such that there is a 
  vertex $w$ in $N(v) \setminus\{u\}$. By Lemma 
  \ref{lemma:fold}, we can remove $w$ from $G$,
  and by induction $G\setminus \{w\}$ can be 
  reduced to a graph without adjacent edges.
\end{proof}


\begin{proposition}\label{prop:addEdge}
  Let $G$ be a graph, $v,w$ distinct non-adjacent
  vertices of $G$, and $G'$ the graph $G$ 
  extended with an edge between $v$ and $w$. 
  If $\mathrm{Ind}(G')$ is $k$--connected and 
  $\mathrm{Ind}(G'\setminus (N_G(v)\cup N_G(w)))$
  is $(k-1)$--connected, then $\mathrm{Ind}(G)$ 
  is $k$--connected.
\end{proposition}
\begin{proof}
  Recall the gluing lemma 
  \cite[10.3(ii)]{Bj}; if $\Delta_1$ and
  $\Delta_2$ are $k$--connected, and $\Delta_1
  \cap \Delta_2$ is $(k-1)$--connected, then
  $\Delta_1 \cup \Delta_2$ is $k$--connected.
  Let $\Delta_1=\mathrm{Ind}(G')$ and $\Delta_2
  =\mathrm{Ind}(G\setminus(N_G(v)\cup N_G(w)))$.
  The complex $\Delta_1$ is $k$--connected by
  assumption, and $\Delta_2$ is $k$--connected 
  since it is a cone. Using $\Delta_1=\{\sigma
  \in\mathrm{Ind}(G)\mid \{v,w\} \not\subseteq
  \sigma\}$ and $\Delta_2=\{\sigma\in\mathrm{Ind}
  (G)\mid\sigma\cup\{v,w\}\in\mathrm{Ind}(G)\}$, 
  we get that $\Delta_1\cup\Delta_2=\mathrm{Ind}
  (G)$, and $\Delta_1\cap\Delta_2=\mathrm{Ind}
  (G'\setminus (N_G(v)\cup N_G(w)))$ which is 
  $(k-1)$--connected. The result follows from 
  the gluing lemma.
\end{proof}


\begin{definition}
  A set of maximal simplices from a simplicial 
  complex $\Delta$ are \emph{generating 
  simplices} if the removal of their interiors
  makes $\Delta$ contractible. Analogously, 
  if $\mathcal{G}$ is a set of maximal faces of 
  $\Delta$ such that $\Delta\setminus\mathcal{G}$
  is contractible, then $\mathcal{G}$ are 
  \emph{generating faces} of $\Delta$.
\end{definition}
  Note that if $\mathcal{G}$ are generating faces
  of $\Delta$, then $\Delta\simeq \vee_{\sigma\in
  \mathcal{G}} S^{\mathrm{dim}\:\sigma}$. A 
  shellable complex has generating faces, they
  are exactly the ones glued over their whole 
  boundary when added. It is not hard to find 
  complexes with generating faces that are not 
  shellable, or to find complexes without 
  generating faces. If $\mathcal{G}$ are 
  generating faces of $\Delta'$, $\Delta$ 
  collapses onto $\Delta'$, and all $\sigma\in 
  \mathcal{G}$ are maximal in $\Delta$, then 
  $\mathcal{G}$ are generating faces of $\Delta$.

  To calculate the homotopy type of complexes
  in this section, a suitable subcomplex
  is found and contracted. Of course one can
  smash any contractible subcomplex, but the
  resulting identifications can be ugly. 
  Our main vehicle is this lemma by Bj\"orner.

\begin{lemma} \label{lemma:homBj}
  \cite[Theorem 10.4(ii)]{Bj}
  Let $\Delta=\Delta_0\cup \Delta_1\cup \cdots 
  \cup \Delta_n$ be a simplicial complex with 
  subcomplexes $\Delta_i$. If $\Delta_i\cap
  \Delta_j\subseteq \Delta_0$ for all $1\leq i<
  j\leq n$, and $\Delta_i$ is contractible for 
  all $0\leq i \leq n$, then
  \[\Delta \simeq \bigvee_{1\leq i\leq n}
  \mathrm{susp} (\Delta_0\cap\Delta_i).  \]
\end{lemma}
  Lemma \ref{lemma:homBj} is a special case of
  \cite[Theorem 2.1]{Bj2}, which can be used to
  generalize both Theorem \ref{theorem:Ind1} and 
  \ref{theorem:Ind2}.
  However, the amount of technicalities do
  not match the increased number of applications 
  at this point.

  Two degenerate cases working well with
  Lemma \ref{lemma:homBj} are that
  $\mathrm{susp}(\emptyset) = S^0$, and that the 
  wedge of nothing is a point.
\begin{theorem}\label{theorem:Ind1}
  If all vertices in the neighborhood of $u\in 
  V(G)$ are adjacent then
  \[ \mathrm{Ind}(G)\simeq \bigvee_{v\in N(u)}
  \mathrm{susp}\: \mathrm{Ind}( G\setminus (
  N(u) \cup N(v) )),\]
  and the union of
  \[ \bigcup_{{v\in N(u) \atop G\setminus (
  N(u) \cup N(v) ) = \emptyset}} \{\{v\}\} \quad
  \quad \mathrm{and} \bigcup_{{v\in N(v) \atop 
  G\setminus (N(u) \cup N(v)) \neq \emptyset}}
  \left\{ \{v\} \cup\sigma \mid \sigma\in
  \mathcal{G}_v\right\} \]
  are generating faces of $\mathrm{Ind}(G)$,
  if $\mathcal{G}_v$ are generating
  faces of $\mathrm{Ind}(G\setminus (
  N(u) \cup N(v)))$.
\end{theorem}
\begin{proof}
  We prove this by smashing the star of $u$
  with Lemma \ref{lemma:homBj} in two ways.
  First let $\Delta_u=
  \mathrm{st}_{\mathrm{Ind}(G)}(u)$, and
  $\Delta_v=\mathrm{st}_{\mathrm{Ind}(G)}(v)$
  for all $v\in N(u)$. Clearly their union
  is $\mathrm{Ind}(G)$, and they
  are all contractible. If $v$ and $w$ are
  different vertices in $N(u)$, then
  $\Delta_v\cap\Delta_w =
  \mathrm{st}_{\mathrm{Ind}(G)}(v) \cap
  \mathrm{st}_{\mathrm{Ind}(G)}(w)
  = v\ast\mathrm{lk}_{\mathrm{Ind}(G)}(v) \cap
  w\ast\mathrm{lk}_{\mathrm{Ind}(G)}(w)
  = \mathrm{lk}_{\mathrm{Ind}(G)}(v) \cap
  \mathrm{lk}_{\mathrm{Ind}(G)}(w)$ since
  $v$ and $w$ are adjacent. Using that
  $N(u)\setminus\{v\} \subseteq 
  N(v)\setminus\{u\}$ and $N(u)\setminus\{w\} 
  \subseteq N(w)\setminus\{u\}$, we get that
  $ \mathrm{lk}_{\mathrm{Ind}(G)}(v) \cap
  \mathrm{lk}_{\mathrm{Ind}(G)}(w) \subseteq
  \mathrm{st}_{\mathrm{Ind}(G)}(u) =\Delta_u$.
  The conditions of Lemma \ref{lemma:homBj} are
  satisfied, and therefore
  \[\begin{array}{rcl}
  \mathrm{Ind}(G) & \simeq & \bigvee_{v\in N(u)}
  \mathrm{susp} (\Delta_u \cap \Delta_v) \\
  &=& \bigvee_{v\in N(u)} \mathrm{susp}\:
  \mathrm{Ind}(G\setminus (N(u) \cup N(v))).\\
  \end{array}\]
  Now to the second part of the theorem.
  Let $V=\{v\in N(u) \mid G\setminus (N(u) \cup 
  N(v)) \neq \emptyset\}$.
  Assume that $\mathcal{G}_v$ are generating
  faces of $\mathrm{Ind}(G\setminus (
  N(u) \cup N(v)))$ for $v\in V$. Let 
  $\mathcal{G}$ be the union of
  \[ \bigcup_{v\in N(u)\setminus V} \{\{v\}\} 
  \quad\textrm{and}\quad \bigcup_{v\in V}
  \left\{ \{v\} \cup\sigma \mid \sigma\in
  \mathcal{G}_v\right\}. \]

  To begin with, we need to prove that each
  element of $\mathcal{G}$ is a maximal face
  of $\mathrm{Ind}(G)$. Let $\sigma$ be a
  maximal face of $\mathrm{Ind}(G\setminus (
  N(u) \cup N(v)))$. Then $\{v\}\cup\sigma$
  is maximal in $\mathrm{Ind}(G\setminus (
  N(u) \setminus \{v\} ))$, and also in
  $\mathrm{Ind}(G)$, since all vertices in
  $N(u) \setminus \{v\}$ are adjacent to $v$.
  If $G\setminus (N(u) \cup N(v))=\emptyset$
  for a $v\in N(u)$, then $N(v)=V(G)\setminus
  \{v\}$, and $\{v\}$ is maximal in
  $\mathrm{Ind}(G)$.
 
  To conclude that $\mathcal{G}$ are 
  generating faces of $\mathrm{Ind}(G)$
  we also need to prove that $\Delta'=
  \Delta\setminus \mathcal{G}$ is 
  contractible. For each $v\in V$ let
  $\Delta_v'= v \ast (\mathrm{Ind}(G\setminus(
  N(u) \cup N(v)))\setminus \mathcal{G}_v )$. 
  It is not hard to see that $\Delta'$ is
  the union of $\Delta_u$ and $\cup_{v\in V}
  \Delta_v'$. All $\Delta_v'$ are contractible,
  and for different $v,w\in V$, $\Delta_v'
  \cap \Delta_w' \subseteq
  \Delta_v \cap \Delta_w \subseteq \Delta_u$.
  By Lemma \ref{lemma:homBj}
  \[\begin{array}{rcl}
  \Delta' & \simeq &
  \bigvee_{v\in V} \mathrm{susp}(
  \Delta_u\cap \Delta_v')\\
  & = & \bigvee_{v\in V} \mathrm{susp}(
  \mathrm{st}_{\mathrm{Ind}(G)}(u)
  \cap v \ast (\mathrm{Ind}(G\setminus (
  N(u) \cup N(v)))\setminus \mathcal{G}_v ) )\\
  & = & \bigvee_{v\in V} \mathrm{susp}(
  \mathrm{Ind}(G\setminus (N(u) \cup N(v)))
  \setminus \mathcal{G}_v  )\\
  & \simeq & \bigvee_{v\in V} 
  \mathrm{susp}(\mathrm{point})\\
  & \simeq & \bigvee_{v\in V} \mathrm{point} \\
  & \simeq & \mathrm{point} \\
  \end{array}\]
  Thus $\mathcal{G}$ are generating faces of 
  $\mathrm{Ind}(G)$.
\end{proof}

  In \cite{K1} the complex $\mathcal{L}
  _n^k$ was defined as the independence complex
  of the graph with vertex set $\{1,2,\ldots 
  n\}$, and two vertices $i<j$ are adjacent if 
  $j-i<k$. The homotopy type $\mathcal{L}_n^k$,
  as well as its generating faces, was 
  calculated 
  for $k=2$, and for $k>2$ stated as an open 
  question. For convenience define 
  $\mathcal{L}_n^k=\emptyset$ if $n\leq 0$.

\begin{corollary}\label{corollary:L}
  For all $k\geq 2$ and $n>1$
  \[\mathcal{L}_n^k \simeq \bigvee
  _{1\leq i< \min(k,n)}\mathrm{susp}
  \left(\mathcal{L}_{n-k-i}^k \right) \]
\end{corollary}
\begin{proof}
  The neighborhood of $1$ is $\{2,3,\ldots 
  ,\mathrm{min}(k,n)\}$, and any two different 
  vertices of it are adjacent.
  \[\begin{array}{rcl}
  \mathcal{L}_n^k & \simeq & \bigvee_{i\in N(1)}
  \mathrm{susp} ( \mathrm{Ind}( G\setminus (
  N(1) \cup N(i) ))) \\
  & = & \bigvee_{i=2}^{\mathrm{min}(k,n)}
  \mathrm{susp} ( \mathrm{Ind}(
  G[\{ j  \mid  k+i \leq j \leq n  \}]))\\
  & \simeq & \bigvee_{i=2}^{\mathrm{min}(k,n)}
  \mathrm{susp} ( \mathrm{Ind}(
  G[\{ j  \mid  1 \leq j \leq n-k-i+1  \}]))\\
  & = &  \bigvee
  _{i=2}^{\mathrm{min}(k,n)} \mathrm{susp}
  \left(\mathcal{L}_{n-k-i+1}^k \right)\\
  & = &  \bigvee
  _{1\leq i< \min(k,n)}\mathrm{susp}
  \left(\mathcal{L}_{n-k-i}^k \right)
  \end{array}\]
\end{proof}
\begin{example}\label{example:L3}
  The generating faces of $\mathcal{L}_n^3$
  produced by recursive use of Theorem 
  \ref{theorem:Ind1}, with $u=1$, are
  \[\begin{array}{r|lr|lr|l}
  n & \mathrm{g.f.} & 
  4 & \{2\},\{3\} &
  8 & \{2,6\},\{2,7\},\{3,7\},\{3,8\} \\
  1 & \emptyset &
  5 & \{3\} &
  9 & \{2,7\},\{3,7\},\{3,8\} \\
  2 & \{2\} &
  6 & \{2,6\} &
  10& \{3,8\},\{2,6,10\} \\
  3 & \{2\},\{3\} &
  7 & \{2,6\},\{2,7\},\{3,7\} &
  11& \{2,6,10\},\{2,6,11\},\{2,7,11\},
      \{3,7,11\}\\
  \end{array}\]
\end{example}


\begin{corollary}\label{corollary:addEdge}
  Let $G$ be a graph with three distinct vertices
  $u,v$ and $w$, such that $N(u)=\{v,w\}$, and
  $\{v,w\}\not\in E(G).$ If $\mathrm{Ind}(G
  \setminus (N(u)\cup N(v)))$ and $\mathrm{Ind}
  (G\setminus (N(u)\cup N(w)))$ are 
  $(k-1)$--connected, and $\mathrm{Ind}(G
  \setminus(N(u)\cup N(v)\cup N(w)))$ is 
  $(k-2)$--connected, then $\mathrm{Ind}(G)$ is 
  $k$--connected.
\end{corollary}
\begin{proof}
  Let $G'$ be the graph $G$ extended with an edge
  between $v$ and $w$. By Proposition 
  \ref{prop:addEdge} it suffices to prove that
  $\mathrm{Ind}(G')$ is $k$--connected and
  $\mathrm{Ind}(G'\setminus(N_G(v)\cup N_G(w)))$
  is $(k-1)$--connected.
  The neighborhood of $u$ in $G'$ is a complete
  graph, so by Theorem \ref{theorem:Ind1},
  \[\begin{array}{rcl}
  \mathrm{Ind}(G') &
  \simeq &
  \mathrm{susp}(
  \mathrm{Ind}(G'\setminus (N_G(u)\cup N_G(v))))
  \vee
  \mathrm{susp}(
  \mathrm{Ind}(G'\setminus (N_G(u)\cup N_G(w))))
  \\
  &=&
  \mathrm{susp}(
  \mathrm{Ind}(G\setminus (N_G(u)\cup N_G(v))))
  \vee
  \mathrm{susp}(
  \mathrm{Ind}(G\setminus (N_G(u)\cup N_G(w)))).
  \\
  \end{array}\]
  The suspension of a $(k-1)$--connected complex
  is $k$--connected, and the wedge of
  $k$--connected complexes is $k$--connected,
  thus $\mathrm{Ind}(G')$ is
  $k$--connected.
  The neighborhood of $v$ in $G'\setminus
  (N_G(v)\cup N_G(w))$ is a complete graph, so
  once again by Theorem \ref{theorem:Ind1},
  \[\begin{array}{rcl}
  \mathrm{Ind}(G'\setminus(N_G(v)\cup N_G(w)))
  &\simeq &
  \mathrm{susp}( \mathrm{Ind}((G' \setminus
  (N_G(v)\cup N_G(w)))\setminus \{v,w\}))  \\
  &=&
  \mathrm{susp}( \mathrm{Ind}(G \setminus
  (N_G(u) \cup N_G(v)\cup N_G(w)))). \\
  \end{array}\]
  The suspension of a $(k-2)$--connected complex
  is $(k-1)$--connected, hence $\mathrm{Ind}(
  G'\setminus (N_G(u)\cup N_G(v)))$ is
  $(k-1)$--connnected.
\end{proof}


  The previous theorem can be used when we find
  a complete subgraph of $G$ with a vertex 
  without neighbours outside the subgraph. 
  Removing the condition of the special vertex 
  forces other conditions. 
\begin{theorem}\label{theorem:Ind2}
  Let $K$ be a subset of $V(G)$ such that $G[K]$ 
  is a complete graph, and $\mathcal{G}$ are 
  generating faces of $\mathrm{Ind}(G\setminus 
  K)$, such that for each $k\in K$ and $\sigma
  \in \mathcal{G}$, one of vertices in $\sigma$ 
  is adjacent to $k$. Then
  \[\mathrm{Ind}(G)\simeq \mathrm{Ind}(G
  \setminus K) \vee\bigvee_{k\in K} \mathrm{susp}
  \: \mathrm{Ind}(G\setminus (K\cup N(k))).\]
  Let $K'=\{k\in K  \mid  G\setminus (K\cup N(k))
  \neq \emptyset \}$. If $\mathcal{G}_k$ are
  generating faces of $\mathrm{Ind}(G\setminus
  (K\cup N(k)))$ for each $k\in K'$, then the 
  union of \[\mathcal{G},\quad
  \bigcup_{k\in K \setminus K'} \{\{k\}\},
  \quad\text{and} \quad
  \bigcup_{k\in K}\{\{k\}\cup\sigma
  \vert\sigma\in\mathcal{G}_k\}\]
  are generating faces of $\mathrm{Ind}(G)$.
\end{theorem}
\begin{proof}
  The proof is in the same spirit as that of
  Theorem \ref{theorem:Ind1}. The subcomplex
  $\Delta_0=\mathrm{Ind}(G\setminus K) \setminus
  \mathcal{G}$ will be contracted. Let
  $\Delta_k=\mathrm{Ind}(G\setminus N(k))$ for
  each $k\in K$, and $\Delta_\tau =
  \{\sigma  \mid  \emptyset \neq \sigma \subseteq 
  \tau \}$ for each $\tau\in \mathcal{G}$.

  If $\sigma \in \mathrm{Ind}(G)$ does not
  contain any vertex from $K$, then $\sigma \in 
  \Delta_\sigma$ if $\sigma\in \mathcal{G}$, and
  $\sigma \in \Delta_0$ if $\sigma\not \in 
  \mathcal{G}$.
  If $\sigma \in \mathrm{Ind}(G)$
  and $k\in\sigma$ for a  $k\in K$, then 
  $\sigma \in \Delta_k$. Hence the union of
  these subcomplexes is $\Delta$.

  Now we check that the required intersections
  are subcomplexes of $\Delta_0$. Note that if 
  $\sigma\in \mathcal{G}$ and $k\in K$, then 
  $\sigma\not \in \Delta_k$ since by assumption 
  there is a vertex in $\sigma$ adjacent to $k$. 
  If $k_1$ and $k_2$ are two different elements 
  of $K$, then $\Delta_{k_1}\cap\Delta_{k_2}
  \subseteq \mathrm{Ind}(G\setminus K)$ since
  $k_1$ and $k_2$ are adjacent. Since $\sigma
  \not\in \Delta_{k_1}$ for any $\sigma\in
  \mathcal{G}$, $\Delta_{k_1}\cap \Delta_{k_2} 
  \subseteq \mathrm{Ind} (G\setminus K)\setminus
  \mathcal{G}=\Delta_0$. If $k\in K$ and $\sigma
  \in \mathcal{G}$, then $\Delta_k\cap
  \Delta_\sigma \subseteq \mathrm{Ind}(G\setminus
  K)$ since $\sigma \in\mathrm{Ind}(G\setminus K)
  $, and $\Delta_k\cap\Delta_\sigma \subseteq 
  \mathrm{Ind}(G\setminus K) \setminus\mathcal{G}
  =\Delta_0$ since $\tau \not\in \Delta_0$ for
  all $\tau \in \mathcal{G}$. If $\sigma_1$ and 
  $\sigma_2$ are different elements of $\mathcal{
  G}$, then $\Delta_{\sigma_1}\cap \Delta_{
  \sigma_2} \subseteq \Delta_0$ since 
  $\sigma_1\not\in \Delta_{\sigma_2}$ and
  $\sigma_2\not\in \Delta_{\sigma_1}$.

  By Lemma \ref{lemma:homBj}
  \[\mathrm{Ind}(G)\simeq\left(\bigvee_{k\in K}
  \mathrm{susp}(\Delta_0 \cap \Delta_k)\right)
  \bigvee\left(\bigvee_{\sigma \in \mathcal{G}}
  \mathrm{susp}(\Delta_0\cap \Delta_\sigma)
  \right).\]
  For $\sigma\in\mathcal{G}$, $\Delta_0 \cap
  \Delta_\sigma = \Delta_\sigma \setminus
  \{ \sigma \} \simeq S^{\mathrm{dim}\sigma -1}$.
  Hence 
  \[\bigvee_{\sigma \in \mathcal{G}}
  \mathrm{susp}(\Delta_0 \cap \Delta_\sigma)
  \simeq \bigvee_{\sigma \in \mathcal{G}}
  \mathrm{susp}\: S^{\:\mathrm{dim}\sigma -1}
  \simeq \bigvee_{\sigma \in \mathcal{G}}
  S^{\:\mathrm{dim}\sigma}
  \simeq \mathrm{Ind}(G\setminus K).\]

  For all $k\in K$, $\Delta_0\cap \Delta_k
  = ( \mathrm{Ind}(G\setminus K)\setminus 
    \mathcal{G}) \cap \mathrm{Ind}
    (G\setminus N(k)) =\mathrm{Ind}(G\setminus K)
    \cap \mathrm{Ind} (G\setminus N(k))
    = \mathrm{Ind} (G\setminus (K\cup N(k)))$
  since for any $\sigma\in \mathcal{G}$ there
  is a $v\in\sigma$ adjacent to $k$, which 
  implies that $\sigma\not\in \mathrm{Ind}
    (G\setminus N(k))$.  
  Inserting this in the conclusion
  of the lemma proves the first part of the
  theorem.

  Now the second part. Let $\mathcal{H}_1=
  \mathcal{G}$, $\mathcal{H}_2= \cup_{k\in K
  \setminus K'}\{\{k\}\}$, and $\mathcal{H}_3=
  \cup_{k\in K}\{\{k\}\cup \sigma \vert
  \sigma \in \mathcal{G}_k\}$. To show that
  $\mathcal{H}=\mathcal{H}_1 \cup \mathcal{H}_2
  \cup \mathcal{H}_3$ are generating faces of
  $\mathrm{Ind}(G)$, we need that all $\sigma\in
  \mathcal{H}$ are maximal faces of $\mathrm{Ind}
  (G)$, and that $\mathrm{Ind}(G)\setminus 
  \mathcal{H}$ is contractible.

  If $k\in K\setminus K'$, then $K\cup N(k)=V(G)
  $. The neighborhood of $k$ is $V(G)\setminus
  \{k\}$ since $K\setminus \{k\}\subseteq N(k)$,
  and $k$ is an isolated point in $\mathrm{Ind}
  (G)$. Thus all elements of $\mathcal{H}_2$ are
  maximal faces of $\mathrm{Ind}(G)$.
  If $\sigma\in\mathcal{H}_1=\mathcal{G}$, then
  $\sigma$ is a maximal face of $\mathrm{Ind}
  (G\setminus K)$. For each vertex $k\in K$
  there is a vertex of $\sigma$ adjacent to it
  by assumption, so no vertex of $K$ can be
  added to $\sigma$. Hence $\sigma$ is also a 
  maximal
  face of $\mathrm{Ind}(G)$.
  Therefore, all elements
  of $\mathcal{H}_1$ are maximal faces of
  $\mathrm{Ind}(G)$.
  If $k\in K'$ and $\sigma \in \mathcal{G}_k$,
  then $\sigma$ is a maximal face of $\mathrm{
  Ind}(G\setminus(K\cup N(k)))=\mathrm{lk}_{
  \mathrm{Ind}(G)}(k)$, so $\{k\}\cup\sigma$
  is a maximal face of $\mathrm{Ind}(G)$. All
  elements of $\mathcal{H}_3$ are therefore
  maximal faces.

  Let $\Delta_0$ be as before, that is $\mathrm{
  Ind}(G\setminus K)\setminus \mathcal{G}$.
  For $k\in K'$ let $\Delta_k'=k\ast(\mathrm{Ind}
  (G\setminus(K\cup N(k)))\setminus \mathcal{G}_k
  )$. We will use Lemma \ref{lemma:homBj} with
  $\Delta_0$ and $\Delta_k'$ for $k\in K$,
  which are all contractible.
  First we show 
  that $\mathrm{Ind}(G)\setminus \mathcal{H}=
  \Delta_0 \cup (\cup_{k\in K'} \Delta_k')$.
  Clearly, $\mathrm{Ind}(G)\setminus \mathcal{H}
  \supseteq
  \Delta_0 \cup (\cup_{k\in K'} \Delta_k')$.
  If $\sigma\in \mathrm{Ind}(G)\setminus 
  \mathcal{H}$ and no vertex of $\sigma$ is in
  $K'$, then $\sigma\in\Delta_0$. If $\sigma\in 
  \mathrm{Ind}(G)\setminus \mathcal{H}$ and 
  $k\in\sigma$, where $k\in K'$, then
  $\sigma\in\Delta_k'$. If $k_1,k_2\in K'$ are
  different, then $\Delta_{k_1}'\cap\Delta_{k_2}'
  =k_1\ast(\mathrm{Ind}(G\setminus(K\cup N(k_1)))
  \setminus \mathcal{G}_{k_1})\cap
  k_2\ast(\mathrm{Ind}(G\setminus(K\cup N(k_2)))
  \setminus \mathcal{G}_{k_2})=
  (\mathrm{Ind}(G\setminus(K\cup N(k_1)))
  \setminus \mathcal{G}_{k_1})\cap
  (\mathrm{Ind}(G\setminus(K\cup N(k_2)))
  \setminus \mathcal{G}_{k_2})\subseteq\Delta_0$.

  By Lemma \ref{lemma:homBj}
  \[\begin{array}{rcl}
  \mathrm{Ind}(G)\setminus\mathcal{H} & \simeq&
  \bigvee_{k\in K'} \textrm{susp}\left(
  \Delta_0 \cap \Delta_k'
  \right)\\
  & = & 
  \bigvee_{k\in K'} \textrm{susp}\left(
  \left(
  \mathrm{Ind}(G\setminus K)\setminus \mathcal{G}
  \right) \cap 
  \left(
  k\ast(\mathrm{Ind}
  (G\setminus(K\cup N(k)))\setminus \mathcal{G}_k
  )
  \right)
  \right)\\
  & = & 
  \bigvee_{k\in K'} \textrm{susp}\left(
  \left(
  \mathrm{Ind}(G\setminus K)\setminus \mathcal{G}
  \right) \cap 
  \left(
  \mathrm{Ind}
  (G\setminus(K\cup N(k)))\setminus \mathcal{G}_k
  \right)
  \right)\\
  & = & 
  \bigvee_{k\in K'} \textrm{susp}\left(
  \mathrm{Ind}(G\setminus K)
  \cap 
  \left(
  \mathrm{Ind}
  (G\setminus(K\cup N(k)))\setminus \mathcal{G}_k
  \right)
  \right)\\
  & = & 
  \bigvee_{k\in K'} \textrm{susp}\left(
  \mathrm{Ind}
  (G\setminus(K\cup N(k)))\setminus \mathcal{G}_k
  \right)\\
  & \simeq & 
  \bigvee_{k\in K'} \textrm{susp}\left(
  \text{point}
  \right)\\  
  & \simeq & 
  \bigvee_{k\in K'} \text{point} \\
  & \simeq & 
  \text{point.} \\
  \end{array}\]
  The equalities need clarification. The first
  one is by definition. The second one follows
  from the fact that $k\not\in \mathrm{Ind}
  (G\setminus K)
  \setminus \mathcal{G}$. Pick a generating
  face $\sigma \in \mathcal{G}$. By assumption,
  there is a $v\in\sigma$ for every $k\in K'$,
  such that $v$ and $k$ are adjacent, that is
  $v\in N(k)$. Thus $\sigma \not \in
  \mathrm{Ind}(G\setminus(K\cup N(k)))\setminus 
  \mathcal{G}_k$, which gives the next equality.
  The final one follows from
  $\mathrm{Ind}(G\setminus K) \supseteq
  \mathrm{Ind}(G\setminus(K\cup N(k)))
  \supseteq \mathrm{Ind}(G\setminus(K\cup N(k)))
  \setminus \mathcal{G}_k$.
\end{proof}
  A relative of $\mathcal{L}_n^k$ is its
  cycle version $\mathcal{C}_n^k$. It is the 
  independence complex of a graph with 
  vertex set $\{1,2,\ldots n\}$, and two 
  vertices $i<j$ are adjacent if $j-i<k$
  or $(n+i)-j<k$. The case $k=2$ was computed
  in \cite{K1}, and used by Babson and 
  Kozlov \cite{BK} in the proof of the 
  Lov\'asz conjecture. The 
  $\mathbb{Z}_2$--homotopy types of 
  $\mathcal{L}_n^2$ and  $\mathcal{C}_n^2$
  were studied by \v{Z}ivaljevi\'c \cite{Z}.
\begin{example}
  Using the generating faces of 
  $\mathcal{L}^3_n$ listed in Example 
  \ref{example:L3}, and Theorem
  \ref{theorem:Ind2}, with $K=\{1,2\}$, we get 
  these generating faces for $\mathcal{C}^3_n$:
  \[\begin{array}{r|l}
  n & \text{generating faces}\\
  \hline
  8  & \{1,5\},\{1,6\},\{2,6\},\{2,7\},\{4,8\} \\
  9  & \{1,5\},\{1,8\},\{2,6\},\{2,9\},
     \{4,8\},\{4,9\},\{5,9\}\\ 
  13 &\{1,5,9\},\{1,5,10\},\{1,6,10\},\{1,6,11\},
     \{2,6,10\},\{2,6,11\},\\
  &    \{2,7,11\},\{2,7,12\},
  \{4,8,12\},\{4,8,13\},\{4,9,13\},\{5,9,13\}\\
\end{array}\]
  Thus $\mathcal{C}^3_8$ is a wedge of five 
  $S^1$, $\mathcal{C}^3_9$ is a wedge of six 
  $S^1$, and $\mathcal{C}^3_{13}$ is a wedge of 
  twelve $S^2$.
\end{example}


  This theorem is molded after Theorem 1.1 in
  \cite{BLVZ}.  
\begin{theorem}\label{Ind3}
  If $G$ is a graph with $n$ vertices and
  maximal degree $d$, then $\mathrm{Ind}(G)$ is
  $\lfloor (n-1)/2d -1 \rfloor$--connected.
\end{theorem}
\begin{proof}
  The proof is by induction on $n$. If $1\leq
  n \leq 2d$ then $\lfloor (n-1)/2d -1 \rfloor
  =-1$ and $\mathrm{Ind}(G)$ is 
  $(-1)$--connected since it is nonempty.

  Recall \cite[Theorem 10.6(ii)]{Bj}: If $\Delta$
  is a simplicial complex and $\{\Delta_i\}
  _{i\in I}$ is a family of subcomplexes such
  that $\Delta=\cup_{i\in I}\Delta_i$, and
  every nonempty intersection $\Delta_{i_1}
  \cap \Delta_{i_2} \cap \cdots \cap \Delta_
  {i_t}$ is $(k-t+1)$--connected, then $\Delta$
  is $k$--connected if and only if the nerve
  $\mathcal{N}(\Delta_i)$ is $k$--connected.

  If $n>2d$ define $\Delta_v=\mathrm{Ind}(G
  \setminus N(v))$ for each $v\in V(G)$. Clearly
  $\Delta=\cup_{v\in V(G)} \Delta_v$.
  The complex $\Delta_v$ is a cone with apex $v$
  and thus $\lfloor (n-1)/2d -1 
  \rfloor$--connected. Let $T$ be a subset of
  $V(G)$ with $t\geq 2$ elements. There are
  at most $d$ vertices in a neighborhood and
  \[ \bigcap_{v\in T} \Delta_v =
  \bigcap_{v\in T} \mathrm{Ind}(G\setminus
  N(v)) = \mathrm{Ind}\left( G\setminus
  \bigcup_{v\in T} N(v)\right),\]
  so $G\setminus \cup_{v\in T} N(v)$ has at
  least $n-td$ vertices and $\mathrm{Ind}(
  G\setminus \cup_{v\in T} N(v))$ is
  $\lfloor (n-td-1)/2d -1 \rfloor$--connected by
  induction. For $t\geq 2$
  \[ \left \lfloor \frac{n-1}{2d}-1 \right 
  \rfloor-t+1 = \left \lfloor \frac{n-td-1}{2d}-
  \frac{t}{2} \right \rfloor \leq \left \lfloor 
  \frac{n-td-1}{2d}-1 \right \rfloor,\]
  thus $\cap_{v\in T}\Delta_v$ is $(\lfloor (n-1)
  /2d-1\rfloor -t+1)$--connected as required.
  We need to show that the nerve is
  $\lfloor (n-1)/2d-1\rfloor$--connected, and
  it will follow from that the intersection
  of $\lfloor (n-1)/2d-1\rfloor+2$ arbitrary
  $\Delta_v$ is nonempty. Indeed, if $T$ is a
  subset of $V(G)$ with $\lfloor (n-1)/2d-1
  \rfloor+2$ elements, then $G\setminus
  \cup_{v\in T} N(v)$ has at least $n-d(\lfloor 
  (n-1)/2d-1\rfloor+2)$ vertices, and
  \[n-d\left( \left \lfloor \frac{ n-1}{2d}-1
  \right \rfloor +2\right) \geq
  n-d\left(  \frac{ n-1}{2d}-1 +2\right)
  = \frac{n-2d}{2}+\frac{1}{2}>\frac{1}{2},\]
  so $\cap_{v\in T}\Delta_v = \mathrm{Ind}(
  G\setminus \cup_{v\in T} N(v))$ is nonempty.
  The conditions of \cite[Theorem 10.6(ii)]{Bj}
  are checked and thus $\mathrm{Ind}(G)$ is
  $\lfloor (n-1)/2d-1 \rfloor$--connected.
\end{proof}
  The independence complex of $m$ disjoint
  complete bipartite graphs $K_{d,d}$ can be 
  collapsed onto the independence complex of 
  $m$ disjoint edges using 
  Lemma \ref{lemma:fold}. That complex is 
  homotopy
  equivalent with an $m-1$ dimensional sphere.
  The $m$ disjoint $K_{d,d}$ have $2md$ vertices
  and maximal degree $d$, thus by 
  Theorem \ref{Ind3} the independence complex is 
  $(m-2)$--connected, which is optimal.

\section{Anti-Rips complexes}
  A natural interpretation of $\mathcal{L}_n^k$
  is as the complex on $\{1,2,\dots n\}\subset 
  \mathbb{R}$, with two different points $p$
  and $q$ in the same simplex if, and only if, 
  $|p-q|> k-1$. Most independence complexes in
  literature can be placed in a metric space,
  which give rise to this definition.
\begin{definition} 
  Let $P$ be a subset of a metric space 
  with distance function $d$, and $r\geq 0$. 
  The \emph{anti-Rips complex} 
  $\mathrm{AR}_r(P)$ have vertex set $P$, and 
  two different points $p$ and $q$ of $P$, 
  can only be in the same simplex if $d(p,q)>r$.
\end{definition} 
  Equivalently, $\mathrm{AR}_r(P)=\mathrm{Ind}
  (G)$, where $G$ is the graph with vertex set 
  $P$, and two different points $p$ and $q$ are 
  adjacent if $d(p,q)\leq r$. Notice that moving 
  $r$ from $0$ to $\infty$ creates a family
  of complexes which is ordered by inclusion,
  and its limits are the simplex on $P$, and $P$
  as disjoint points.

  Why name it anti-Rips complexes? Substituting
  $d(p,q)>r$ with $d(p,q)\leq r$ defines
  Rips complexes. According to Hausmann 
  \cite{H}
  Lefschetz called them Vietoris complexes, but
  the notation changed with Rips' reintroduction
  of them in the study of hyperbolic groups.
  A contemporary application of Rips complexes 
  is in the approximation of homotopy type of
  point-cloud data, see for example Carlsson
  and da Silva \cite{CdS}. Corollary 
  \ref{corollary:L} can now be generalized.
\begin{proposition}
  If $P$ is a finite subset of $\mathbb{R}$
  and $m=\min(P)$ then
  \[\mathrm{AR}_r(P) \simeq \bigvee_
  {p\in P \atop m<p\leq m+r}
  \mathrm{susp}(\mathrm{AR}_r(\{
  q\in P  \mid  q>p+r\}))\]
\end{proposition}
\begin{proof}
  Let $u=m$ in Theorem \ref{theorem:Ind1}.
  Use that the neighborhood of a point $p$ in
  the graph corresponding to $\mathrm{AR}_r(P)$,
  is the set of points $q$ in $P$, such that
  $d(p,q)>r$.
\end{proof}

  If $P\subset \mathbb{Z}^2$ have $n$ vertices,
  then $\mathrm{AR}_1(P)$ is $\lfloor (n-9)/8
  \rfloor$--connected by Theorem \ref{Ind3}
  since there are at most 4 points within
  distance one from a point in $\mathbb{Z}^2$.
  Using the geometry of the plane,
  $\lfloor (n-9)/8 \rfloor$ can be improved
  to $\lfloor (n-9)/6 \rfloor$.
\begin{proposition}
  If $P\subset \mathbb{Z}^2$ have $n$ vertices,
  then $\mathrm{AR}_1(P)$ is $\lfloor (n-9)/6
  \rfloor$--connected.
\end{proposition}

\begin{proof}
  The proof is by induction on $n$. If $1\leq n
  \leq 8$ then $\mathrm{AR}_1(P)$ is $\lfloor
  (n-9)/6 \rfloor$--connected since it is 
  $(-1)$--connected, and $\lfloor (n-9)/6\rfloor
  \leq -1$.
  Now assume that $n>8$. Pick a $u\in P$ such
  that the sum of its $x$ and $y$ coordinates
  is maximal among the points in $P$. It is no
  restriction to assume that $u=(0,0)$ since
  the proposition is translation invariant.
  Define $N:P\rightarrow 2^P$ by $N(p)=\{q\in P
  \mid d(p,q)\leq 1, p\neq q\}.$ Let $v=(0,-1)$
  and $w=(-1,0)$. Depending on $N(u)\subseteq
  \{v,w\}$ we have four cases.

  If $N(u)=\emptyset$ then $\mathrm{AR}_1(P)$ is
  a cone and in particular $\lfloor (n-9)/6
  \rfloor$--connected.

  If $N(u)=\{v\}$ then $\mathrm{AR}_1(P)\simeq
  \mathrm{susp}(\mathrm{AR}_1(P\setminus
  (N(u)\cup N(v))))$ by Theorem 
  \ref{theorem:Ind1}.
  The complex $\mathrm{AR}_1(P\setminus (N(u)\cup
  N(v)))$ is $\lfloor ((n-5)-9)/6
  \rfloor$--connected by induction since
  $N(u)\cup N(v)\subseteq \{u,v,(1,-1),(0,-2),
  (-1,-1)\}.$ The suspension increases the
  connectivity by one, and 
  $\lfloor ((n-5)-9)/6  \rfloor+1 \geq \lfloor 
  (n-9)/6 \rfloor$, thus
  $\mathrm{AR}_1(P)$ is $\lfloor (n-9)/6
  \rfloor$--connected. The case $N(u)=\{w\}$ is
  analogous.
 
  The final case is $N(u)=\{v,w\}$. By induction
  $\mathrm{AR}_1(P \setminus (N(u) \cup N(v)
  \cup N(w)))$ is $(\lfloor (n-9)/6
  \rfloor-2)$--connected since $N(u) \cup N(v)
  \cup N(w)\subseteq \{u,v,w, (-1,1), (-2,0),$
  $(-1,-1), (1,-1), (0,-2) \}$ and 
  $\lfloor ((n-8)-9)/6 \rfloor \geq 
  \lfloor (n-9)/6 \rfloor-2.$ 
  The complex $\mathrm{AR}_1(P \setminus 
  (N(u) \cup N(v)))$ is $(\lfloor
  (n-9)/6 \rfloor-1)$--connected by induction
  since $N(u) \cup N(v) \subseteq \{u,v,w,
  (-1,-1), (1,-1), (0,-2) \}.$ Analogously
  $\mathrm{AR}_1(P \setminus (N(u) \cup N(w)))$ 
  is $(\lfloor (n-9)/6 \rfloor-1)$--connected.
  By Corollary \ref{corollary:addEdge}, 
  $\mathrm{AR}_1(P)$ is $\lfloor (n-9)/6 
  \rfloor$--connected.
\end{proof}


\section{Open questions}
  We conclude with some open questions.
\begin{question}
  One approach to bound the connectivity
  of a simplicial complex is to chop it up
  in pieces for which the connectivity can be 
  calculated easily, and then use the Nerve Lemma
  (cf. \cite{BLVZ,E} and Theorem \ref{Ind3}).
  Suitable subcomplexes for complexes of
  directed trees are those with the same roots
  and Theorem \ref{theorem:homEqDAG} shows
  their connectedness. However their 
  intersections are in general cumbersome.
  Can this class of subcomplexes be
  adapted to prove a nontrivial
  bound for the connectivity?
\end{question}
\begin{question}
  Theorem \ref{theorem:Ind2} puts conditions
  on the generating faces which in general
  are hard to verify, and maybe even not
  possible to achieve by choosing the
  generating faces correctly. In practice, when
  all generating faces are ``far from a certain
  vertex'', that vertex can often be 
  collapsed away, or discrete Morse theory 
  \cite{F} can be used. Can this be formalized
  to a method for removing vertices not in 
  generating faces?
\end{question}
\begin{question}
  Kozlov \cite{K2} found explicit homotopy
  equivalences between $\mathcal{L}^2_{2n}$
  and complexes described by Shapiro and
  Welker \cite{SW}. Is it possible to
  generalize to $\mathcal{L}^k_n$,
  or even to anti-Rips complexes
  on $\mathbb{R}$?
\end{question}
\begin{question}
  The homotopy type of $\mathcal{C}^k_n$ in 
  general is still unsolved. A larger class
  to investigate is the anti-Rips complex of
  a finite subset of a circle.
\end{question}

\end{document}